\documentclass[a4paper]{article}
\usepackage{graphicx}
\usepackage{amsmath}
\usepackage{amsfonts}
\usepackage{amsthm}
\usepackage{amssymb}
\usepackage{subfigure}
\usepackage{multirow}
\usepackage{rotating}
\usepackage{fancybox}
\usepackage{boxedminipage}
\usepackage{version}
\usepackage[usenames]{color}
\usepackage[english,greek]{babel}
\usepackage[iso-8859-7]{inputenc}
\newcommand{\ds}{\displaystyle}
\newcommand{\en}{\selectlanguage{english}}

\newcommand{\NN}{\mathbb{N}}
\newcommand{\RR}{\mathbb{R}}
\newcommand{\ZZ}{\mathbb{Z}}

\newcommand{\FF}{\mathbb{F}}

\newcommand{\QQ}{\mathbb{Q}}
\newcommand{\be}{\begin{equation}}
\newcommand{\ee}{\end{equation}}

\newtheorem{thm}{Theorem}

\theoremstyle{definition}
\newtheorem{rmk}{Remark}
\newtheorem{dfn}{Definition}

\newtheorem{cor}{Corollary}

\newtheorem{alg}{Algorithm}

\begin{document} \en
\title{On the generalization of the Costas property in the continuum}
\author{Konstantinos Drakakis\footnote{The author holds a Diploma in Electrical and Computer Engineering from NTUA, Athens, Greece, and a Ph.D. in Applied and Computational Mathematics from Princeton University, NJ, USA. He was a scholar of the Lilian Boudouris Foundation.}, Scott Rickard
\\Electronic and Electrical Engineering\\University College Dublin\\ \& \\ Claude Shannon Institute\footnote{www.shannoninstitute.ie}\\Ireland}
\maketitle

\abstract{We extend the definition of the Costas property to functions in the continuum, namely on intervals of the reals or the rationals, and argue that such functions can be used in the same applications as discrete Costas arrays. We construct Costas bijections in the real continuum within the class of piecewise continuously differentiable functions, but our attempts to construct a fractal-like Costas bijection there are successful only under slight but necessary deviations from the usual arithmetic laws. Furthermore, we are able, contingent on the validity of Artin's conjecture, to set up a limiting process according to which sequences of Welch Costas arrays converge to smooth Costas bijections over the reals. The situation over the rationals is different: there, we propose an algorithm of great generality and flexibility for the construction of a Costas fractal bijection. Its success, though, relies heavily on the enumerability of the rationals, and therefore it cannot be generalized over the reals in an obvious way.}

\section{Introduction}
Costas arrays \cite{C} have been an active topic of research for more than 40 years now; however, after 1984, when 2 algebraic construction methods for Costas arrays were published (the Welch and the Golomb method \cite{G}), still the only ones available today, there has been effectively no progress at all in the construction of new Costas arrays, with the obvious exception of brute force searches. Recent research on Costas arrays tends to focus on the discovery of new properties \cite{D2,D3,SVM}, hoping that they will either furnish some lead for a new construction method, or prove that such a method does not exist, and thus overcome the current virtual stalemate in the core problems of the field.

In line with this effort, it is likely that research on Costas arrays would benefit by the extension of the definition of the Costas property in the continuum, for 2 reasons: on the one hand, this might open the door to assistance from the entire arsenal of analysis, as was the case with the successful generalization of the factorial in terms of the Gamma function; on the other hand, the recent advances in the subject of the Instantaneous Frequency of a signal \cite{HSLWSZTL} make it possible to design signals with continuously varying frequencies instead of piecewise constant frequencies, such as the usual discrete Costas arrays model, and there might be benefits in doing so. And besides, such objects certainly have an intrinsic pure mathematical merit for study.

In this work, we propose a suitable extension of the definition of the Costas property in the continuum (which we take here to mean the real and rational numbers), and we explain how the existing discrete Costas permutations can be used to generate continuum Costas permutations. Note that, in accordance with common practice in recent literature, we will be using the terms ``Costas permutation'' and ``Costas array'' interchangeably.

\section{Basics}
We reproduce below the definition of a Costas function/permutation \cite{D}:

\begin{dfn}\label{bas} Let $[n]:=\{0,\ldots,n-1\},\ n\in\NN$ and consider a bijection $f:[n]\rightarrow [n]$; $f$ is a Costas permutation iff the multiset $\{(i-j,f(i)-f(j)): 0\leq j<i< n\}$ is actually a set, namely all of its elements are distinct.
\end{dfn}

These permutations are extremely useful because they give rise to binary signals with an optimal autocorrelation pattern:

\begin{dfn}\label{dac} Let $f:[n]\rightarrow [n]$, $n\in\NN^*$, be a Costas permutation, and let $F:\ZZ^2\rightarrow [2]$, the corresponding binary signal of $f$, satisfy $F(i,f(i))=1,\ i\in[n]$, and $F=0$ everywhere else. The autocorrelation of $f$ is:
\[A_F(u,v)=\sum_{i,j\in\ZZ} F(u+i,v+j)F(i,j),\ (u,v)\in\ZZ^2\]
\end{dfn}
The following result is just a restatement of the Costas property:

\begin{thm}
Let $f:[n]\rightarrow [n]$, $n\in\NN^*$, be a permutation, and let $F$ be its corresponding binary signal; then, $0\leq A_f(u,v)<2,\ \forall u,v\in\ZZ^2-\{0,0\}$ iff $f$ has the Costas property.
\end{thm}

We have already mentioned the Welch construction method for Costas arrays. As we will refer to it several times below, we offer its definition for the sake of completeness:

\begin{thm}[Welch construction $W_1(p,g,c)$] Let $p$ be a prime, let $g$ be a primitive root of the finite field $\FF(p)$ of $p$ elements, and let $c\in[p-1]$ be a constant; then, the function $f:[p-1]+1\rightarrow [p-1]+1$ where $\ds f(i)=g^{i-1+c}\mod p,\ i\in[p-1]$ is a bijection with the Costas property.

\end{thm}

\section{Costas bijections in the real continuum}

From now on, until Section \ref{crat}, we will be using the term ``continuum'' in the sense of ``real continuum'', unless explicitly stated otherwise.

\subsection{Definitions and simple results}\label{defs}

In our extension of Definition \ref{bas} in the continuum we will replace $[n]$ by $[0,1]$, but otherwise the definition remains the same:

\begin{dfn} Consider a bijection $f:[0,1]\rightarrow [0,1]$; $f$ \emph{is a Costas permutation} iff the multiset $\{(x-y,f(x)-f(y)): 0\leq y<x\leq 1\}$ is actually a set, namely all of its elements are distinct.
\end{dfn}

\begin{rmk}
The choice of the interval $[0,1]$ is by no means restrictive: it can be seen immediately that for any pair $a,b\in\RR$, $a<b$ there exists a linear monotonic mapping $h$ mapping $[0,1]$ bijectively on $[a,b]$, specifically $h(x)=a+x(b-a),\ 0\leq x\leq 1$, and $f$ has the Costas property on $[a,b]$ iff $h^{-1}\circ f\circ h$ has the Costas property on $[0,1]$.
\end{rmk}

Yet again, we can give an alternative but equivalent definition of the Costas property in terms of autocorrelation:

\begin{dfn}
Consider a bijection $f:[0,1]\rightarrow [0,1]$, and let $F:\RR^2\rightarrow \{0,1,\infty\}$ be its corresponding quasi-binary signal (that is, binary whenever finite), so that $F(x,f(x))=1,\ x\in[0,1]$, and $F=0$ otherwise. The autocorrelation of $f$ is:
\[A_f(u,v)=\int_0^1 \int_0^1 \delta(F(x+u,y+v)-F(x,y))dxdy,\ (u,v)\in\RR^2\]
\end{dfn}

\begin{rmk}
Notice that this autocorrelation, just like its discrete counterpart in Definition \ref{dac}, takes integer values whenever finite, as it counts the number of zeros in the argument of the Dirac $\delta$-function.
\end{rmk}

Once more, then, the following result is just a restatement of the Costas property:

\begin{thm}
Consider a bijection $f:[0,1]\rightarrow [0,1]$, and let $F$ be its corresponding quasi-binary signal; then, $f$ has the Costas property iff $0\leq A_f(u,v)<2,\ \forall (u,v)\in\RR^2-\{0,0\}$.
\end{thm}

\subsection{Applications}

Continuum Costas bijections\footnote{We have to resort to the use of the uncommon word ``continuum'' in the role of an adjective here instead of the perhaps more appealing intuitively ``continuous'': the term ``continuum function'' accurately describes a function defined on an interval, or on something non-finite and dense at any rate, whereas the term ``continuous function'' has an already established different meaning in mathematics.} can find applications in the same situations their discrete counterparts do \cite{C}. For example, consider a RADAR system whose operation relies on a usual Costas waveform. In practical terms, this means that the waveform it transmits is of the form:
\[w(t)=A\cos\left(2\pi\left(\sum_{k=0}^{n-1} \frac{s(k)+1}{n}f \mathbf{1}_{\left[\frac{k}{n}T,\frac{k+1}{n}T\right)}(t)\right)t\right),\ s\text{ a Costas permutation of order $n$},\ t\in[0,T)\]
which is a different way to express that, for $n\in\NN^*$, $\ds w(t)=A\cos\left(2\pi \frac{s(k)+1}{n}ft\right),\ t\in\left[\frac{k}{n}T,\frac{k+1}{n}T\right),\ k\in[n]$.

Alternatively, we could have used a continuum Costas permutation $s$ on $[0,1]$. Let us consider the waveform:
\[w(t)=A\cos\left(2\pi f \int_0^t s(u)du+2\pi f_0 t\right),\ t\in[0,T)\]
Bedrosian's theorem \cite{B, HSLWSZTL} on instantaneous frequency asserts that the instantaneous frequency of $w$ is
\[ \frac{1}{2\pi}\left(2\pi f \int_0^t s(u)du+2\pi f_0 t\right)'=s(t)f+f_0,\]
as long as $\hat{w}(0)=0$; this condition can be satisfied, at least approximately, through an appropriate choice of $f_0$.

\subsection{Link between continuum and discrete Costas permutations}

How do the 2 definitions compare? The expression for the discrete waveform is clearly a special case of the continuum expression, and this can be seen if we write $\ds S(t)=\sum_{k=0}^{n-1} \frac{s(k)+1}{n}f \mathbf{1}_{\left[\frac{k}{n}T,\frac{k+1}{n}T\right)}(t)$, where $s$ is a Costas array of order $n$ and $S$ is a continuum permutation (but obviously not Costas). The verification of the Costas property through the autocorrelation in the discrete case is also a subprocess of the verification in the continuum case: we just need to take care that horizontal and vertical displacements of the copies of the functions in the autocorrelation formula are integral multiples of $\ds \frac{T}{n}$ and $\ds \frac{f}{n}$, respectively.

As $S$ is a piecewise constant function, one might be tempted to attempt to formulate a definition for (at least a class of) continuum Costas permutations in terms of Costas arrays, as limits of sequences of Costas arrays, just like measurable functions are approximated by sequences of piecewise constant functions: a Costas array $s_n$ of order $n$ can be mapped on a piecewise constant function $S_n$, just as we did above, and, letting $n\rightarrow \infty$, we can hopefully obtain a continuum Costas permutation $S$. This limit would probably be highly discontinuous, of a fractal nature perhaps, as Costas arrays are highly erratic and patternless.

The problem with the plan of action suggested above is that we seem to have no good understanding yet of sequences of Costas arrays across different orders that follow a clear pattern, so that we can successfully describe how the limit of such a sequence would look like; a notable exception is the example we give below in Section \ref{lm}. Nevertheless, the idea of seeking continuum Costas permutations among fractals seems, in principle, promising in itself and worthwhile investigating. But first, let us focus on the case of smooth functions.

\section{Construction of smooth continuum Costas permutations}\label{smperm}

The whole idea of the existence of smooth functions with the Costas property may sound outright irrational at first, and any investigation futile: after all, there can hardly be any object more irregular and discontinuous that Costas arrays. Nonetheless, the continuum is dense in itself, while finite discrete sets are not, and this makes a big difference, as we are about to see: for example, the function $f(x)=x^2$ has no chance of being a permutation on any discrete set other than $[2]\cup\{\infty\}$, while it is a permutation on both $[0,1]$ and $[1,+\infty]$, as it effectively makes some areas of the intervals ``denser'' and some ``sparser'' (consider, for instance, that the images under $f$ of all points in $[0,\sqrt{2}/2]$ get ``crammed'' in the smaller interval $[0,0.5]$). In the continuum we can create Costas permutations by causing ``elastic deformations'', by ``changing the density'' of points in an interval, whereas such techniques are inapplicable on discrete sets.

Let us begin by seeking functions with the Costas property that are reasonably smooth; for example, let us confine ourselves to special categories of almost everywhere differentiable bijections.

\begin{dfn}
Let $f:[0,1]\rightarrow [0,1]$ be a bijection;
\begin{itemize}
  \item $f$ will be \emph{piecewise continuously differentiable} iff there exists $n\in\NN^*\cup \{\infty\}$ and a sequence of intervals $\{I_i\}_{i=1}^n$ so that, for each $i=1,\ldots,n$, $f$ is continuously differentiable in $\ds \overset{\circ}{I_i}$ ($n=\infty$ is used as a convention to denote a countable infinity of intervals);
  \item if, in addition to being piecewise continuously differentiable, for each $i=1,\ldots,n$ $\ds f'|\overset{\circ}{I_i}$ is strictly monotonic, $f$ will be called \emph{piecewise strictly monotonic piecewise continuously differentiable};
  \item if, in addition to being piecewise continuously differentiable, $f$ satisfies the property that, for all sequences of points $\{x_i\}_{i=1}^n$ such that $\ds x_i\in \overset{\circ}{I_i},\ i=1,\ldots,n$, it is true that the sequences $\{f'(x_i)\}_{i=1}^n$ are either all strictly increasing or all strictly decreasing, $f$ will be called \emph{overall strictly monotonic piecewise continuously differentiable};
  \item $f$ may combine all 3 features above, in which case it will be called \emph{overall and piecewise strictly monotonic piecewise continuously differentiable}.
\end{itemize}
\end{dfn}

\begin{thm}
Let $f:[0,1]\rightarrow [0,1]$ be an overall and piecewise strictly monotonic continuously differentiable bijection. Then, $f$ has the Costas property on $[0,1]$.
\end{thm}

\begin{proof}
Let us choose 4 points in $[0,1]$, say $x$, $y$, $x+d$ and $y+d$ so that $y< x$ and $d\geq 0$; these may actually be 3 points if $x=y+d$. We need to show that
\[f(x)-f(x+d)=f(y)-f(y+d)\Rightarrow d=0\]
Exactly one of the 2 pairs of intervals $[x,y],[x+d,y+d]$ or $[x,x+d], [y,y+d]$ consists of intervals with disjoint interiors. Without loss of generality, assume it is the second pair, then the Newton-Leibnitz Theorem implies that
\[f(x+d)-f(x)=\int_{x}^{x+d}f'(u)du,\ f(y+d)-f(y)=\int_{y}^{y+d}f'(u)du\]
Now, if $f$ is overall and piecewise strictly monotonic continuously differentiable, it is always the case that either $\forall u,v:\ u\in(x,x+d),\ v\in(y,y+d) f'(u)< f'(v)$ or $\forall u,v:\ u\in(x,x+d),\ v\in(y,y+d) f'(u)> f'(v)$, so that $f(x)-f(x+d)\neq f(y)-f(y+d)$ unless $d=0$.
\end{proof}

\begin{thm}
Let $f:[0,1]\rightarrow [0,1]$ be a piecewise continuously differentiable bijection; if $f'$ is not injective, $f$ does not have the Costas property.
\end{thm}

\begin{proof}
We distinguish the following cases:
\begin{itemize}
  \item $f'$ is constant on an interval, say $f'\equiv c\in\RR$ , or, equivalently, $f$ is linear on that interval: it follows there exist 4 points $x$, $y$, $x+d$ and $y+d$ with $y< x$ and $d> 0$ so that $\ds \frac{f(x+d)-f(x)}{d}=\frac{f(y+d)-f(y)}{d}=c$, hence the Costas property is violated.
  \item Assume that $f'$ is never constant on an interval. Then, either there exist $i_1,i_2$ so that $|f'(I_{i_1})\cap f'(I_{i_2})|>0$, namely it fails to be overall strictly monotonic, or there exists an $i$ for which $f'|I_i$ is not monotonic. In either case, there exist 2 points $x_1,x_2\in(0,1)$, so that $x_1<x_2$ and $f'(x_1)=f'(x_2)$. We distinguish 2 subcases:
  \begin{itemize}
    \item Neither of the points is an inflection point, that is both points lie in regions of the domain where $f$ is either convex or concave; these regions are necessarily different, or the derivative could not possibly be equal at these points. This implies that there exist real numbers $\epsilon_1,\epsilon_2>0$ so that, if 2 parallels are drawn to the tangent at each of the points $x_1$ and $x_2$, at the side of the tangents where the function graph lies, and whose distances from the tangents are less than $\epsilon_1$ and $\epsilon_2$, respectively, they each intersect the function graph at 2 points, say $x_{11}<x_{12}$ and $x_{21}<x_{22}$. Clearly both $x_{11}-x_{12}$ and $x_{21}-x_{22}$ go to 0 as the parallels move closer to the tangents, whence $f(x_{11})-f(x_{12})$ and $f(x_{21})-f(x_{22})$ also go to 0; moreover, if $\epsilon_1$ and $\epsilon_2$ are sufficiently small, $(x_{11},x_{12})\cap (x_{21},x_{22})=\emptyset$, and each of $(x_{11},x_{12}), (x_{21},x_{22})$ falls entirely within one of the intervals $\{I_i\},\ i=1,ldots,n$. Hence, we can choose a pair of parallels so that $\ds \frac{f(x_{11})-f(x_{12})}{x_{11}-x_{12}}= \frac{f(x_{21})-f(x_{22})}{x_{21}-x_{22}}$ and $x_{11}-x_{12}=x_{21}-x_{22}$. This violates the Costas property.
    \item At least one of the points is an inflection point, say $x_1$, so there is a $\delta$ so that $x\in(x_1-\delta,x_1+\delta)-\{x_1\}\Rightarrow f'(x)<f'(x_1)$ and $(x_1-\delta,x_1+\delta)$ falls within one of the intervals $\{I_i\},\ i=1,ldots,n$, say $I_k$. As $f'$ is continuous within $I_k$, and is not constant in any interval, there exist $u_1\in(x_1-\delta,x_1)$, $u_2\in(x_1,x_1+\delta$ so that neither is an inflection point and that $f'(u_1)=f'(u_2)$. We are now back to the case above.
  \end{itemize}
\end{itemize}
\end{proof}

Note that the derivative of a continuously differentiable bijection must keep the same sign throughout its domain, or else the bijection would have an extremum and would not be a bijection. Further, in the case of a continuously differentiable bijection, overall and piecewise strict monotonicity are identical, hence strict monotonicity implies injectivity. Therefore, in this special case, the following holds:

\begin{cor}\
\begin{itemize}
  \item Let $f:[0,1]\rightarrow [0,1]$ be a bijection continuously differentiable in $(0,1)$; then, $f$ has the Costas property iff $f'$ is strictly monotonic.
  \item A continuously differentiable bijection on $f:[0,1]\rightarrow [0,1]$ with the Costas property must be strictly monotonic.
\end{itemize}
\end{cor}

\begin{rmk}
The issue of the continuity of the derivative of a function is quite esoteric. When a function is differentiable in an open interval, its derivative is not necessarily continuous. However, it is ``almost'' continuous, in the sense that, for any value between 2 values the derivative actually assumes at 2 points, there is a point between the 2 aforementioned points where the derivative assumes the chosen value. This property is known as Darboux continuity in the literature \cite{BC}. Working with piecewise continuously differentiable functions, we ``float over'' this technical point.
\end{rmk}

Let us now see some examples of continuously differentiable bijections with the Costas property as well as some rules to produce new ones from known ones:

\begin{cor}\label{exmp}
The following continuously differentiable bijections  $f:[0,1]\rightarrow [0,1]$ have the Costas property on $[0,1]$:
\begin{itemize}
  \item $f(x)=x^a$, $a\in\RR_+$, $a\neq 0,1$;
  \item $\ds f(x)=\frac{a^x-1}{a-1},\ a\in\RR^*_+ -\{1\}$;
  \item $\ds f(x)=\sin\left(\frac{\pi}{2} x\right)$;
\end{itemize}
Further, if $f,g:[0,1]\rightarrow [0,1]$ are continuously differentiable bijections and have the Costas property on $[0,1]$, the following functions also do:
\begin{itemize}
  \item $1-f$;
  \item $af+bg$, $a,b\in\RR_+$, $a+b=1$, if $f,g$ are both strictly increasing or all strictly decreasing, and so are $f',g'$;
  \item $f\circ g$, if $f',g'$ are strictly monotonic of the same type and $g$ is strictly increasing;
  \item $fg$, if $f,g,f',g'$ are all strictly increasing or all strictly decreasing.
\end{itemize}
\end{cor}

\begin{proof}
Observe that $\ds \left(\frac{a^x-1}{a-1}\right)'=\ln(a)\frac{a^x}{a-1}$ is strictly increasing for $a>1$ and strictly decreasing for $a<1$, $(x^a)'=ax^{a-1}>0$ is strictly increasing when $a>1$ and strictly decreasing when $0<a<1$, and $\ds \left(\sin\left(\frac{\pi}{2} x\right)\right)'=\frac{\pi}{2}\cos\left(\frac{\pi}{2} x\right)$ is strictly decreasing. Moreover, all of these functions are bijections, hence, they have the Costas property.

Further,
\begin{itemize}
  \item $(1-f)'=-f'$ is strictly monotonic iff $f'$ is, although of the opposite type, and $1-f$ is a bijection on $[0,1]$, so it also has the Costas property.
  \item $(af+bg)'=af'+bg'$ is strictly monotonic if $f',g'$ are both strictly monotonic of the same type, and $af+bg$ is strictly monotonic too, hence a bijection, if $f,g$ are both strictly monotonic of the same type.
  \item $f\circ g$ is clearly a bijection if both $f$ and $g$ are, and $(f\circ g)' = g'f'\circ g$ is strictly increasing (decreasing) if both $f',g'$ are strictly increasing (decreasing) and $g$ is strictly increasing.
  \item $fg$ is strictly increasing (decreasing), hence a bijection, if $f,g$ are both strictly increasing (decreasing), while $(fg)'=fg'+f'g$ is strictly increasing (decreasing) if $f,g,f',g'$ are all strictly increasing (decreasing).
\end{itemize}
\end{proof}

We have now offered a quite extensive description of the class of piecewise continuously differentiable bijections on $[0,1]$ with the Costas property, and an exact characterization of the continuously differentiable bijections with the Costas property. What about discontinuous bijections, though? By interpreting discontinuity in the most extreme way, we are led back to the idea of fractals.

\section{Costas fractals}\label{cofrac}

In what follows, we establish a connection between discrete and continuum Costas permutations: we use discrete Costas permutations to build continuum ones through a process of multiscale rearrangement of subintervals of $[0,1]$; in other words, we build a ``Costas fractal''. At this moment, however, we are unable to prove the correctness of our construction below under the usual laws of arithmetic: we will need the equivalent of ``xor'' addition (and subtraction), namely addition without carry, in representations over an arbitrary basis.

We will need first of all the slightly stronger definition given below:

\begin{dfn}Consider a bijection $f:[n]\rightarrow [n]$; $f$ is a \emph{modulo Costas permutation} iff the multiset $\{(i-j,f(i)-f(j)\mod(n+1)): 0\leq j<i< n\}$ is actually a set, namely all of its elements are distinct.
\end{dfn}

\begin{rmk} Note that both the Golomb and the Welch constructions actually lead to modulo Costas permutations \cite{D,G}.
\end{rmk}

\begin{dfn}\label{no}
Let the numbers $x,y\in[0,1]$ be expanded over basis $n\in\NN^*$: $\ds x=\sum_{i=1}^\infty x_i n^{-i}$, $\ds y=\sum_{i=1}^\infty y_i n^{-i}$, where $\forall i\in\NN^*, x_i,y_i\in[n]$. Then, we define the ``no carry'' addition and subtraction as:
\[x\oplus y=\sum_{i=1}^\infty \frac{(x_i+y_i)\mod n}{n^i},\ x\ominus y=\sum_{i=1}^\infty \frac{(x_i-y_i)\mod n}{n^i}\]
\end{dfn}

\begin{thm}\label{fr}
Let $n\in\NN$ and let $f_i:[n]\rightarrow [n],\ i\in\NN^*$ be a sequence of (not necessarily distinct) modulo Costas permutations. Define a function $F:[0,1]\rightarrow [0,1]$ by the following formula:
\[F\left(\sum_{i=1}^\infty a_i n^{-i}\right)=\sum_{i=1}^\infty f_i(a_i)n^{-i}\]
where $\forall i\in\NN^*, a_i\in[n]$, and so that there exists no $N\in\NN^*: a_i=n-1$ for $i\geq N$, unless $N=1$. Then, $F$ has the Costas property, when subtraction is interpreted as in Definition \ref{no}.
\end{thm}

\begin{rmk}
The explicit exclusion of sequences $\{a_i\}_{i=1}^\infty$ so that $\exists N\in\NN^*: a_i=n-1$ for $i\geq N$ is necessary in order to ensure that every number in $[0,1)$ can be expressed over base $n$ in a unique way, otherwise some numbers can have 2 different expansions: a familiar example over base 10 would be that $0.5=0.5000\ldots=0.4999\ldots$. However, we still need to represent $\ds 1=\sum_{i=1}^\infty \frac{n-1}{n^i}$, hence the exception for $N=1$. \end{rmk}

\begin{proof}
Select 4 points in $[0,1]$, say $x$, $y$, $x+d$ and $y+d$ so that $y< x$ and $d\geq 0$; notice that these can actually be 3 equidistant points if $y+d=x$. We need to test whether $F(x)\ominus F(x+ d)=F(y)\ominus F(y+ d)$ necessarily implies $d=0$.

Let the interval $[0,1]$ be divided into $n$ subintervals, $\ds \left\{I_{1;i}=\left[\frac{i}{n}, \frac{i+1}{n}\right):i\in[n-1]\right\}\bigcup \left\{I_{1;n-1}=\left[\frac{n-1}{n}, 1\right] \right\}$, so that $\forall i\in[n], F(I_{1;i})=I_{1;f(i)}$. We distinguish the following cases:

\begin{enumerate}
  \item \label{main} $y+d\neq x$ and the 4 chosen points all lie in different subintervals: then, we can write $\ds F(x)=\frac{s_1}{n}+\epsilon_1$, $\ds F(y)=\frac{s_2}{n}+\epsilon_2$, $\ds F(x+d)=\frac{s_3}{n}+\epsilon_3$, and $\ds F(y+d)=\frac{s_4}{n}+\epsilon_4$, with $s_i\in[n]$, $\ds \epsilon_i<\frac{1}{n},\ i=1,2,3,4$. It follows that $\ds F(x)\ominus F(x+d)=\frac{(s_1-s_3)\mod n}{n}+(\epsilon_1 \ominus \epsilon_3)$, and $\ds F(y)\ominus F(y+d)=\frac{(s_2-s_4)\mod n}{n}+(\epsilon_2 \ominus \epsilon_4)$, where, if we assume $d>0$, $(s_1-s_3)\mod n\neq (s_2-s_4) \mod n $, by the modulo Costas property of $f_1$, while $\ds |(\epsilon_1 \ominus \epsilon_3)\ominus (\epsilon_2 \ominus \epsilon_4)|<\frac{1}{n}$. Hence, $F(x)\ominus F(x+d)\neq F(y)\ominus F(y+d)$ and the proof is complete for this case.
  \item $y+d= x$ and the 3 chosen points all lie in different subintervals: then we can repeat verbatim the previous argument with 3 instead of 4 points.
  \item $y+d\neq x$ and one pair of the 4 chosen points lie in the same subinterval, while the remaining pair lie in different subintervals: then, without loss of generality, assume that $x$ and $x+d$ lie in the same subinterval. In terms of the previous argument, $(s_1-s_3)\mod n=0\neq (s_2-s_4)\mod n$ and the proof follows again.
  \item $y+d= x$ and the 3 chosen points lie in 2 different subintervals: then, exactly 2 points lie in the same subinterval, and, without loss of generality, assume they are $y$ and $y+d=x$. In terms of the previous argument, $s_4=s_1$, $(s_1-s_3)\mod n\neq (s_2-s_1)\mod n=0$ and the proof follows again.
  \item\label{bad} Either $y+d\neq x$ and the 4 chosen points lie pairwise in the same subintervals, or $y+d= x$ and the 3 chosen points all lie in the same subinterval: then, assume, without loss of generality, that $x$ and $x+d$ lie in the same subinterval, and so do $y$ and $y+d$. It follows that $(s_1-s_3)\mod n=0= (s_2-s_4)\mod n$ and the argument fails.
\end{enumerate}

In the last case where the argument fails, we need to refine our subinterval division. We already saw the first level of this division. At level $k\in\NN$, we consider the collection of intervals
\begin{multline*}
\left\{I_{k;i_1,\ldots,i_k}=\left[\sum_{j=1}^k\frac{i_j}{n^j}, \sum_{j=1}^{k-1}\frac{i_j}{n^j}+\frac{i_k+1}{n^j}\right):\ i_j\in[n],\ j=1,\ldots,k,\ \exists j:i_j\neq n-1\right\}\bigcup\\
\left\{I_{k;n-1,\ldots,n-1}=\left[1-\frac{1}{n^k},1\right]\right\}
\end{multline*}

With respect to the newly defined levels of subintervals, there are 2 possibilities:
\begin{itemize}
  \item The chosen points fall in a case other than \ref{bad} for the first time in level $k$: then, it must be the case that:
   \begin{multline*} \sum_{j=1}^k \frac{(f_j(x_j+d_j)-f_j(x_j))\mod n}{n^j}=\frac{(f_k(x_k+d_k)-f_k(x_k))\mod n}{n^k}\neq\\ \sum_{j=1}^k \frac{f_j((y_j+d_j)-f_j(y_j))\mod n}{n^j}=\frac{(f_k(y_k+d_k)-f_k(y_k))\mod n }{n^k}\end{multline*}
   due to the modulo Costas property of $f_k$, whence $F(x)\ominus F(x+d)\neq F(y)\ominus F(y+d)$ for $d>0$.
  \item Otherwise, we need to consider the levels beyond level $k$.
\end{itemize}
But the length of the subintervals in level $k$ is $n^{-k}$ which decays to 0 as $k\rightarrow \infty$; therefore, any specific selection of points can remain in case \ref{bad} for a finite number of levels only. This completes the proof.
\end{proof}

It is easy to see where our proof fails under ordinary arithmetic: revisiting case \ref{main}, we would need to show that, under the assumption that $s_1-s_3\neq s_2-s_4$, which holds because $f_1$ is a Costas permutation (we no longer need it to be a modulo Costas permutation),  $\ds \frac{s_1-s_3}{n}+(\epsilon_1 - \epsilon_3)\neq \frac{s_2-s_4}{n}+(\epsilon_2 - \epsilon_4)$ holds. Since $\ds \epsilon_i<\frac{1}{n},\ i=1,2,3,4$, it follows that $\ds |\epsilon_1 - \epsilon_3|,|\epsilon_2 - \epsilon_4|<\frac{1}{n}$ and $\ds |(\epsilon_1- \epsilon_3)- (\epsilon_2 - \epsilon_4)|<\frac{2}{n}$, so that, if $|(s_1-s_3)-(s_2-s_4)|=1$, it may still be the case that $\ds \frac{s_1-s_3}{n}+(\epsilon_1 - \epsilon_3)= \frac{s_2-s_4}{n}+(\epsilon_2 - \epsilon_4)\Leftrightarrow F(x)-F(x+d)=F(y)-F(y+d)$ when $d>0$, and the Costas property fails.

The key feature of the arithmetic proposed in Definition \ref{no} that allowed the proof of Theorem \ref{fr} to complete successfully was that if, at any level of interval subdivision, the 4 chosen points were found to lie into distinct subintervals, the defining inequality of the Costas property would be satisfied for the chosen points. There are alternative arithmetics with this property:

\begin{dfn}\label{np}

Let the numbers $x,y\in[0,1]$ be expanded over basis $n\in\NN^*$: $\ds x=\sum_{i=1}^\infty x_i n^{-i}$, $\ds y=\sum_{i=1}^\infty y_i n^{-i}$, where $\forall i\in\NN^*, x_i,y_i\in[n]$. Then, we define the ``contracted'' subtraction as:
\[x\ominus y=\sum_{i=1}^\infty \frac{x_i-y_i}{n^{2i-1}}\]
\end{dfn}

\begin{thm}\label{fs}
Let $n\in\NN$ and let $f_i:[n]\rightarrow [n],\ i\in\NN^*$ be a sequence of (not necessarily distinct) Costas permutations. Define a function $F:[0,1]\rightarrow [0,1]$ by the following formula:
\[F\left(\sum_{i=1}^\infty a_i n^{-i}\right)=\sum_{i=1}^\infty f_i(a_i)n^{-i}\]
where $\forall i\in\NN^*, a_i\in[n]$, and so that there exists no $N\in\NN^*: a_i=n-1$ for $i\geq N$, unless $N=1$. Then, $F$ has the Costas property, when subtraction is interpreted as in Definition \ref{np}.
\end{thm}

\begin{proof}
This is a verbatim repetition of the proof of Theorem \ref{fr}.
\end{proof}

Is it likely that Theorem \ref{fr} still hold true for ordinary arithmetic despite the fact that our proof does not carry through? At this time we have no reason to believe that it does. It may still be possible to use discrete Costas permutations to generate a Costas fractal in the continuum, but the actual mechanism should most probably be different.

\section{A limiting process}\label{lm}

Assuming Artin's Conjecture holds true \cite{M}, which would be the case if the Generalized Riemann Hypothesis holds true, for any non-square integer $k\in\NN^*$, $k>1$ there exists an infinite sequence of primes, say $\{p_n\}_{n\in\NN^*}$, for which $k$ is a primitive root. We can construct then the sequence of Welch Costas permutations corresponding to the primes of the sequence and the primitive root $k$:
\[f_n: [p_n-1]+1\rightarrow [p_n-1]+1,\ f_n(i)=k^{i-1}\mod p_n,\ i\in[p_n-1]+1,\ n\in\NN^*\]

The key observation is that $\forall m\in\NN^*, \exists N\in\NN^*:\ \forall i=1,\ldots,m, f_n(i)=k^{i-1}$; in particular, $N$ is the smallest integer for which $p_N>k^{m-1}$. In other words, for any fixed number of terms, all functions of the sequence, after skipping a finite number of functions, have these initial terms in common. Define then the \emph{pointwise intermediate limit} of $\{f_n\}_{n\in\NN^*}$ to be as follows: for a fixed $i\in\NN^*$,
\[f(i)=\lim f_n(i)=\lim_{n>N} f_n(i):=\lim k^{i-1},\text{ where $N$ is the smallest integer such that $p_N>k^{i-1}$}\]

Choose now a sequence $\{i_n\}_{n\in\NN^*}$ of integers such that $\ds \lim\frac{i_n-1}{p_n}=x$. We define the limit of $\{f_n\}_{n\in\NN^*}$ evaluated on $\{i_n\}_{n\in\NN^*}$ to be a continuum function on $[0,1]$  as follows:
\[s(x)=\lim \left(f(i_n)\right)^{\frac{1}{p_n}}=k^x,\ x\in[0,1)\]
We can bring the range of $s$ within $[0,1)$ as well after a linear transformation, and create: $\ds S(x)=\frac{k^x-1}{k-1}$; this is the second example function in Corollary \ref{exmp}.

To sum up, in the special case of an infinite sequence of Welch Costas permutations generated by a common primitive root $k$, we were able to carry out a limiting process and construct a continuum Costas permutation, using the property that all the members of this sequence (except possibly some of the first ones) have a common beginning. The limit we obtained, however, is a smooth function and not a fractal, as one might expect given the way Welch Costas permutations look like.

\section{Costas bijections in the rational continuum} \label{crat}

The idea of fractals with the Costas property in the (real) continuum was explored above in Section \ref{cofrac}, where we saw that their implementation required special considerations. We return to this issue here, but this time in the context of the rationals $Q=\QQ\cap [0,1]$: in many ways the rationals stand midway between the integers and the reals, in the sense that they form a dense set (like the reals), but still enumerable (like the integers). We are about to see that these 2 properties allow us to make further progress in the subject.

Note that Costas permutations on the rational continuum is a genuinely new problem, and in no way a special case of the constructions in the real continuum; the reason is that the constructions of Section \ref{smperm} do not map bijectively the rationals onto the rationals. For example, $f(x)=x^2$ is not a bijection over $Q$, as $\ds \nexists x\in Q:\ f(x)=\frac{1}{3}$, say.

The relevant definitions of the Costas property on rational bijections closely parallel the ones in Section \ref{defs} (regarding the real continuum) and will not be repeated here.

\subsection{An existence result}

In this section we offer an algorithm of considerable generality for the construction of bijections on $Q$ with the Costas property. Let us begin by reordering the elements of $Q$ as follows: we order firstly by the magnitude of the denominator, and secondly by the magnitude of the numerator (both in an increasing way). Explicitly, first come those rational numbers in $[0,1]$ whose denominator is 1, namely $\ds 0=\frac{0}{1}$ and $\ds 1=\frac{1}{1}$; then, those whose denominator is 2, namely $\ds \frac{1}{2}$; then, those whose denominator is 3, namely $\ds \frac{1}{3}$ and $\ds \frac{2}{3}$ etc. Hence, the sequence looks like this:
\[0,1,\frac{1}{2},\frac{1}{3},\frac{2}{3},\frac{1}{4},\frac{3}{4},\frac{1}{5},\frac{2}{5},\frac{3}{5},\frac{4}{5},\frac{1}{6},\frac{5}{6}\ldots \]
Notice that the numerators are always taken to be relatively prime to the denominators in order to avoid duplicate entries. We denote $Q$ equipped with this particular ordering by $Q_X$, and its elements, in the order dictated by the ordering, by $x_0,x_1,x_2,\ldots$. This ordering has the advantage that each rational is preceded by a finite number of rationals only (in set theoretic terminology, it does not contain any transfinite points). Similarly, we denote by $Q_Y$ the set $Q$ equipped with any arbitrary but fixed ordering without transfinite points, and we denote its elements, in the order dictated by its ordering, by $y_0,y_1,y_2,\ldots$.

Consider now the following algorithm for the construction of a mapping $f:Q\rightarrow Q$:

\begin{alg}\ \label{algq}

\begin{description}
	\item[Initialization] Choose $f(x_0)=y_0$; set $Q'_{Y}\leftarrow Q_Y-\{y_0\}$, $Q'_{X}\leftarrow Q_X-\{x_0\}$, $X\leftarrow \{x_0\}$, $Y\leftarrow \{y_0\}$, and $D\leftarrow \{\}$.
	\item[Find $x$ for $y$:] Set $Q_{X,\text{av}}\leftarrow Q'_X$, $x\leftarrow \inf Q_{X,\text{av}}$, $y\leftarrow \inf Q'_{Y}$; while the (multi)set $\{\text{sgn}(x'-x)(x'-x, f(x')-y): x'\in X\}\cup D$ is actually a multiset, set $Q_{X,\text{av}}\leftarrow Q_{X,\text{av}}-\{x\}$, $x\leftarrow \inf Q_{X,\text{av}}$, and repeat. Set $f(x)=y$, $D\leftarrow \{\text{sgn}(x'-x)(x'-x, f(x')-y): x'\in X\}\cup D$, $Q'_{Y}\leftarrow Q'_{Y}-\{y\}$, $Q'_{X}\leftarrow Q'_{X}-\{x\}$, $X\leftarrow X\cup \{x\}$, $Y\leftarrow Y\cup \{y\}$.
    \item[Find $y$ for $x$:] Set $Q_{Y,\text{av}}\leftarrow Q'_Y$, $y\leftarrow \inf Q_{Y,\text{av}}$, $x\leftarrow \inf Q'_{X}$; while the (multi)set $\{\text{sgn}(x'-x)(x'-x, f(x')-y): x'\in X\}\cup D$ is actually a multiset, set $Q_{Y,\text{av}}\leftarrow Q_{Y,\text{av}}-\{y\}$, $y\leftarrow \inf Q_{Y,\text{av}}$, and repeat. Set $f(x)=y$, $D\leftarrow \{\text{sgn}(x'-x)(x'-x, f(x')-y): x'\in X\}\cup D$, $Q'_{Y}\leftarrow Q'_{Y}-\{y\}$, $Q'_{X}\leftarrow Q'_{X}-\{x\}$, $X\leftarrow X\cup \{x\}$, $Y\leftarrow Y\cup \{y\}$.
\end{description}
The algorithm needs to be supplied with a step sequence before execution begins. For the purposes of the correctness proof the exact step sequence is unimportant (this is yet another degree of freedom of the algorithm), as long as the following rules are observed:
\begin{itemize}
  \item Initialization is run first and only once;
  \item Neither Find $x$ for $y$ nor Find $y$ for $x$ is run infinitely many times in a row.
\end{itemize}
\end{alg}

For example, when $Q_Y=Q_X$ and the steps are run alternatingly, we get $\ds f(0)=0,\ f(1)=1,\ f\left(\frac{1}{2}\right)=\frac{1}{3},\ f\left(\frac{1}{3}\right)=\frac{1}{2},\ f\left(\frac{2}{3}\right)=\frac{2}{3}$ etc.

\begin{thm}\label{thmq}
Algorithm \ref{algq} produces infinitely many bijections $f:Q\rightarrow Q$ with the Costas property.
\end{thm}

\begin{proof}
In order to prove the correctness of Algorithm \ref{algq} above, we need to demonstrate that a) $\forall y\in Q_Y, \exists! x\in Q_X: f(x)=y$, and b) $\forall x\in Q_X, \exists y\in Q_Y: f(x)=y$. To begin with, note that the construction algorithm above guarantees that the constructed $f$ has the Costas property and that every $y\in Q_Y$ appears in the range of $f$ at most once. We only need to show that the algorithm never gets ``stuck'', namely that the two while loops always exit.
\begin{itemize}
  \item For a given $x$, is it possible to assign a value to $f(x)$? In other words, if $A\subset Q'_Y$ is the set of all values $f(x)$ can take without violating the Costas property of $f$, is it true that $A\neq \emptyset$? The answer is in the affirmative, as, intuitively, we can see that the Costas property restrictions impose only a finite number of constraints on $f(x_i)$, while $Q'_Y$ is countably infinite. Rigorously, we have to check 2 conditions:
  \begin{itemize}
    \item Let $A_1\subset Q'_Y$ be the set of possible values for $f(x)$ for which $\text{sgn}(x'-x)(x'-x,f(x')-f(x))=\text{sgn}(x'-x'')(x'-x'',f(x')-f(x''))$ is never true for $x',x''\in X$. We show that $A_1\neq \emptyset$. In fact, consider $\ds \frac{1}{p}$, where $p$ is a prime that does not appear as a factor in the denominator of some $f(x'),\ x'\in X$: choosing $\ds f(x)=\frac{1}{p}$, it follows that $\ds \frac{1}{p}-f(x')$ contains $p$ as a factor in the denominator, while $f(x')-f(x'')$ does not, hence they cannot be equal, and therefore that $\ds \frac{1}{p}\in A_1\neq \emptyset$ as promised. Clearly, there are infinitely many choices for $p$ possible, so $A_1$ contains actually infinitely many elements.
    \item Let $A_2\subset A_1$ be the set of possible values for $f(x)$ for which $\text{sgn}(x'-x)(x'-x,f(x')-f(x))=\text{sgn}(x''-x)(x''-x,f(x'')-f(x))$ is never true for $x',x''\in X$. We show that $A_2\neq \emptyset$. In order for one of these equalities to hold, $x$ must be the midpoint of $x'$ and $x''$, while at the same time $f(x)$ be the midpoint of $f(x')$ and $f(x'')$. Choosing $\ds f(x)=\frac{1}{p}$ where $p$ is as above, and writing $\ds x'=\frac{u_1}{v_1},\ x''=\frac{u_2}{v_2}$, we need to investigate whether the following is possible:
    \[\frac{1}{2}\left(\frac{u_1}{v_1}+\frac{u_2}{v_2}\right)=\frac{1}{p},\ (u_1,v_1)=(u_2,v_2)=1,\ p\not|v_1,v_2.\]
    This implies $p(u_1 v_2+u_2 v_1)=2v_1v_2$, and therefore that $p|2v_1v_2\Rightarrow p|2\Rightarrow p=2$. Hence, $A_2$ does contain all points of the form $\ds \frac{1}{p}$, too, where $p$ does not divide the denominator of some $f(x'),\ x'\in X$ (which are infinitely many), except possibly for $\ds \frac{1}{2}$; in any case $A_2\neq \emptyset$.
  \end{itemize}
  But $A_2=A$, hence $A\neq \emptyset$, a contradiction; therefore, $f(x)$ can assume a value without $f$ losing the Costas property.

  \item For a given $y$, is it possible to find $x\in Q'_X$ so that $f(x)=y$? In other words, if $A\subset Q'_X$ is the set of all values $x$ for which $f(x)$ can be $y$ without violating the Costas property of $f$, is it true that $A\neq \emptyset$? The answer is in the affirmative as well, and the argument is an almost verbatim repetition of the argument above. Rigorously, we have to check 2 conditions:
  \begin{itemize}
    \item Let $A_1\subset Q'_X$ be the set of possible values for $x$ for which $\text{sgn}(x'-x)(x'-x,f(x')-y)=\text{sgn}(x'-x'')(x'-x'',f(x')-f(x''))$ is never true for $x',x''\in X$. We show that $A_1\neq \emptyset$. In fact, consider $\ds \frac{1}{p}$, where $p$ is a prime that does not appear as a factor in the denominator of some $x'\in X$: choosing $\ds x=\frac{1}{p}$, it follows that $\ds \frac{1}{p}-x'$ contains $p$ as a factor in the denominator, while $x'-x''$ does not, hence they cannot be equal, and therefore that $\ds \frac{1}{p}\in A_1\neq \emptyset$ as promised. Clearly, there are infinitely many choices for $p$ possible, so $A_1$ contains actually infinitely many elements.
    \item Let $A_2\subset A_1$ be the set of possible values for $x$ for which $\text{sgn}(x'-x)(x'-x,f(x')-y)=\text{sgn}(x''-x)(x''-x,f(x'')-y)$ is never true for $x',x''\in X$. We show that $A_2\neq \emptyset$. In order for one of these equalities to hold, $x$ must be the midpoint of $x'$ and $x''$, while at the same time $y$ be the midpoint of $f(x')$ and $f(x'')$. Choosing $\ds x=\frac{1}{p}$ where $p$ is as above, and writing $\ds x'=\frac{u_1}{v_1},\ x''=\frac{u_2}{v_2}$, we need to investigate whether the following is possible:
    \[\frac{1}{2}\left(\frac{u_1}{v_1}+\frac{u_2}{v_2}\right)=\frac{1}{p},\ (u_1,v_1)=(u_2,v_2)=1,\ p\not|v_1,v_2.\]
    This implies $p(u_1 v_2+u_2 v_1)=2v_1v_2$, and therefore that $p|2v_1v_2\Rightarrow p|2\Rightarrow p=2$. Hence, $A_2$ does contain all points of the form $\ds \frac{1}{p}$, too, where $p$ does not divide the denominator of some $f(x'),\ x'\in X$ (which are infinitely many), except possibly for $\ds \frac{1}{2}$; in any case $A_2\neq \emptyset$.
  \end{itemize}
  But $A_2=A$, hence $A\neq \emptyset$, a contradiction; therefore, there exists a $x:f(x)=y$ without $f$ losing the Costas property.
\end{itemize}
This completes the proof.
\end{proof}

\begin{rmk}
Intuitively, the mechanism responsible for the flexibility of the algorithm is the opportunity the countable infinity of the rationals offers for ``double deference of all difficulties for a future time'': when faced with the difficulty of assigning a value to $f$ at a given point, we always have infinitely many possibilities, out of which some will work; this in turn creates the difficulty of assigning the values we skipped to some point, but, when faced with this difficulty, we again have infinitely many points waiting for an assignment, out of which some again will work; but in choosing one we once more skip some points, and we need to choose values for them, hence the cycle restarts.

This interplay is precisely what we cannot do with a finite set, hence the contrast between the easiness of the Costas construction over the rationals, as opposed to the intractability of the classical construction of Costas arrays.
\end{rmk}

\begin{rmk}
The above proof makes heavy use of the enumerability of the rationals, and therefore cannot be readily extended to the reals, who lack this property.
\end{rmk}

It may come as a surprise that we can extend the algorithm even further:

\begin{thm}
Algorithm \ref{algq} will produce a bijection $f:Q\rightarrow Q$ with the Costas property even if one of the steps From $x$ to $y$ or From $y$ to $x$ is applied infinitely many times in a row.
\end{thm}

\begin{proof}
Let us consider the case where From $x$ to $y$ is run infinitely many times in a row immediately after Initialization. This causes no loss of generality: the case where From $y$ to $x$ is run infinitely many times in a row immediately after Initialization is completely dual (observe the duality in the proof of Theorem \ref{thmq}), while the more general situation where finitely many alternations between the 2 steps occur before the algorithm ``locks'' in one can be considered to fall within one of the 2 cases we just mentioned, but with a different, more extensive Initialization.

Assume then that we go through $x\in Q_X$ one after another and we try to assign values to $f(x)\in Q_Y$ while retaining the Costas property. The proof of Theorem \ref{thmq} guarantees that we will succeed for all points. What we need to worry about is whether some $y\in Q_Y$ will be left out in the process: in other words, we know that $\forall x\in Q_X,\ \exists y\in Q_Y: f(x)=y$, but we still need to know that $\forall y\in Q_Y,\ \exists! x\in Q_X:\ f(x)=y$.

Assume then that at some step of the algorithm we find that $y\in Q'_Y$ has been skipped, and is the smallest element of $Q_Y$ that has been skipped. Will the algorithm ever ``pick it up''? As before, let us denote by $A\subset Q'_X$ the set of all available $x$ for which we can set $f(x)=y$ without violating the Costas property; we need to show that $A\neq \emptyset$. Because we proceed through $Q_X$ sequentially from the beginning, at the particular step of the algorithm we find ourselves there exists $x_0\in Q_X:\ X=\{x\in Q_X: x\leq x_0\}$ (remember that $\leq$ refers to the ordering of $Q_X$, \emph{not} the usual ordering!).

Consider a $x\in Q'_X$ which is of the form $\ds \frac{1}{p}$, $p$ prime, say $\chi$. As in the proof of Theorem \ref{thmq}, we need to show 2 things:
\begin{itemize}
	\item $\text{sgn}(x'-\chi)(x'-\chi,f(x')-y)=\text{sgn}(x'-x'')(x'-x'',f(x')-f(x''))$ is never true for $x',x''<\chi$. The additional complication here is that at the current step of the algorithm we know the values of $f$ up to $x_0$, but we endeavor to prove a property that holds for $x<\chi$, i.e. involving future values! The way to avoid the complication is to apply our favorite argument on the first coordinate only, disregarding entirely what the values of $f$ are: $x'-x''$ cannot contain $p$ as a factor in its denominator, while $\chi-x'$ does, hence they cannot be equal. It follows that $\chi$ will belong in $A$ as long as it satisfies the second condition we are now about to test, and also that $\chi$ can actually be chosen among infinitely many points.
	\item $\text{sgn}(x'-\chi)(x'-\chi,f(x')-y)=\text{sgn}(x''-\chi)(x''-\chi,f(x'')-y)$ is never true for $x',x''<\chi$. In order to check this we repeat verbatim the proof of Theorem \ref{thmq}: we assume that $\chi$ is the midpoint of some $x'$ and $x''$, and then show this is impossible, unless perhaps $\ds \chi=\frac{1}{2}$. It follows that $\ds \frac{1}{p}$ satisfies this condition too, with the possible exception of when $p=2$. But this still leaves infinitely many points of the form $\ds \frac{1}{p}$, $p$ prime, in $A$, hence in particular $A\neq \emptyset$.
\end{itemize}
This completes the proof.
\end{proof}

\subsection{An explicit construction}

Algorithm \ref{algq} is not exactly constructive; we cannot, for example, readily compute what $\ds f\left(\frac{8}{1025}\right)$ is equal to. We propose here a constructive algorithm for the construction of a Costas permutation on the rationals; the catch is, however, that it only works on a subset of $Q$.

\begin{dfn}
We define the set of \emph{prime rationals} $Q_P$ in $[0,1]$ to be the subset of $Q$ with prime denominators; namely $\ds Q_P=\left\{\frac{i}{p}:\ i\in[p-1],p\text{ prime}\right\}$.
\end{dfn}

\begin{thm}
For each prime $p$, consider a Welch Costas permutation $f_p:[p-1]+1\rightarrow [p-1]+1$ constructed in $\FF(p)$, and consider the set of points $\ds S(p)=\left\{\left(\frac{i}{p},\frac{f_p(i)}{p}\right):\ i\in[p-1]+1\right\}$. The set $\ds S=\bigcup_{p\text{ prime}} S(p)$ is a Costas permutation on $Q_P$.
\end{thm}

\begin{proof}
$S$ is clearly a permutation. We need to show that the distance vectors between all pairs of points are distinct.
\begin{itemize}
  \item Choose 4 points in the same $S(p)$: the Costas property of $f_p$ guarantees the 2 distance vectors they define are distinct.
  \item Choose 2 points in $S(p)$ and 2 points in $S(q)$, $q\neq p$: the first distance vector has coordinates that are fractions over $p$, while the second over $q$, hence they cannot be equal.
  \item Choose 2 points in $S(p)$, a point in $S(q)$, and a point in $S(r)$, where $p,q,r$ are distinct primes: the first distance vector has coordinates that are fractions over $p$, while the second over $qr$, hence they cannot be equal.
  \item Choose a point in $S(p)$, a point in $S(q)$, a point in $S(r)$, and a point in $S(s)$, where $p,q,r,s$ are distinct primes: the first distance vector has coordinates that are fractions over $pq$, while the second over $rs$, hence they cannot be equal.
\end{itemize}
This completes the proof.
\end{proof}

\section{Conclusion}

In this work, we have made 4 main and original contributions to the subject of Costas arrays:
\begin{itemize}
  \item We defined the Costas property on a real continuum function in 2 ways, through distance vectors between points and through the autocorrelation, and we showed that the 2 definitions are equivalent. We also showed that real continuum Costas bijections can be used in the same applications as discrete Costas arrays, by designing signals with the appropriate instantaneous frequency, which has been made possible by the recent advances in the field. Subsequently, we studied similarly the Costas property on rational continuum functions. Essentially, we have now translated the entire framework of Costas arrays in the continuum.
  \item We showed that real continuum Costas bijections exist and we offered some examples; we characterized completely the continuously differentiable Costas bijections in terms of the monotonicity of their derivative, and we also obtained some good results for the case where the bijections are only piecewise continuously differentiable.
  \item We investigated whether it is possible to construct fractal bijections with the Costas property, perhaps by employing discrete Costas arrays as building blocks. We answered that in the affirmative under nonstandard arithmetic laws (where addition and subtraction take place without carry, or where the contribution of the least significant digits of the points to their distance is deemphasized) in the real continuum; under ordinary arithmetic we have no reason to believe that the result still holds true.
  \item We proposed a very general and flexible algorithm for the construction of Costas permutations over the rationals, that is not, however, entirely constructive. We were also able to formulate such a constructive algorithm, but its applicability is limited over a subset of the rationals.
\end{itemize}

Overall, it comes to us as a surprise that it was relatively simple to construct smooth continuum functions with the Costas property, whereas all efforts to create a fractal Costas real bijection were unsuccessful (under ordinary arithmetic). Intuitively, given the irregularity of discrete Costas arrays, we would expect the known construction methods for Costas arrays to generalize in a natural way in the real continuum leading to a fractal; however, a direct recursion, such as our attempt in Section \ref{cofrac}, seems to be inappropriate, unless we change the arithmetic we use. It may still be possible to construct a Costas fractal bijection based on discrete Costas arrays through a different, less obvious mechanism, and we challenge the reader to discover such a mechanism.

\end{document}